\documentclass[11pt,a4paper,leqno,twoside]{amsart}
\usepackage{latexsym,amssymb,amsmath}
\input xy
\xyoption{all}

\def\CC{\mathbb C}

\def\RR{\mathbb R}
\def\HH{\mathbb H}
\def\AA{{\mathbb A}}

\def\OO{\mathbb O}

\def\11{\mathbf 1}
\def\PP{\mathbb P}

\def\e1{\varepsilon_1}
\def\e2{\varepsilon_2}
\def\e3{\varepsilon_3}

\def\P2{{\PP}^2}

\def\00{\underline{0}}
\def\J0{{\cal J}_3(\underline{0})}

\def\PJ0{\PP({\cal J}_3(\underline{0}))}

\def\e{\varepsilon}

\def\AP2{{\AA\PP}^2}
\def\RP2{{\RR\PP}^2}
\def\CP2{{\CC\PP}^2}
\def\HP2{{\HH\PP}^2}
\def\OP2{{\OO\PP}^2}

\newtheorem{theo}{Theorem}[section]

\newtheorem{lemm}[theo]{Lemma}

\theoremstyle{remark}
\newtheorem{rema}[theo]{Remark}
\theoremstyle{definition}

\begin{document}

\title{Lower bounds and infinity criterion for Brauer $p$-dimensions
of finitely-generated field extensions}

\vskip0.6truecm

\author{I.D. Chipchakov}
\par
\medskip
\address{Institute of Mathematics and Informatics\\ Bulgarian
Academy of Sciences\\ Acad. G. Bonchev Str., bl. 8\\ 1113, Sofia,
Bulgaria; chipchak@math.bas.bg} 
\thanks{Partially supported by a project No. RD-08-241/12.03.2013 of 
Shumen University.}

\par
\begin{abstract}
Let $E$ be a field, $p$ a prime number and $F/E$ a
finitely-generated extension of transcendency degree $t$. This paper
shows that if the absolute Galois group $\mathcal{G}_{E}$ is of
nonzero cohomological $p$-dimension cd$_{p}(\mathcal{G}_{E})$, then 
the field $F$ has Brauer $p$-dimension Brd$_{p}(F) \ge t$ except, possibly, in case $p = 2$, the Sylow pro-$2$-subgroups of $\mathcal{G}_{E}$ are of order $2$, and $F$ is a nonreal field. It announces that Brd$_{p}(F)$ is infinite whenever $t \ge 1$ and the absolute Brauer $p$-dimension abrd$_{p}(E)$ is infinite; moreover, 
for each pair $(m, n)$ of integers with $1 \le m \le n$, there exists 
a central division $F$-algebra of exponent $p ^{m}$ and Schur index 
$p ^{n}$.
\end{abstract}

\maketitle

\par
\section{Introduction and index-exponent relations over
finitely-generated field extensions}
\par
\medskip
Let $E$ be a field, Br$(E)$ its Brauer group, $s(E)$ the class of 
finite-dimensional associative central simple $E$-algebras, and 
$d(E)$ the subclass of division algebras $D \in s(E)$. It is known 
that Br$(E)$ is an abelian torsion group (cf. \cite{P}, Sect. 14.4), 
so it decomposes into the direct sum of its $p$-components Br$(E) 
_{p}$, where $p$ runs across the set $\mathbb P$ of prime
numbers. Denote by $[A]$ the equivalence class in Br$(E)$ of any $A 
\in s(E)$. The degree deg$(A)$, the Schur index ind$(A)$, and the
exponent exp$(A)$ (the order of $[A]$ in Br$(E)$) are important
invariants of $A$. Note that deg$(A) = n.{\rm ind}(A)$, and
ind$(A)$ and exp$(A)$ are related as follows (cf. \cite{P}, Sects.
13.4, 14.4 and 15.2):
\par
\medskip
(1.1) exp$(A)$ divides ind$(A)$ and is divisible by every $p \in
\mathbb P$ dividing ind$(A)$. For each $B \in s(E)$ with ind$(B)$
relatively prime to ind$(A)$, ind$(A \otimes _{E} B) = {\rm
ind}(A).{\rm ind}(B)$; in particular, the tensor product $A \otimes
_{E} B$ lies in $d(E)$, provided that $A \in d(E)$ and $B \in d(E)$.
\par
\medskip
As shown by Brauer, (1.1) fully describe the generally valid
restrictions between Schur indices and exponents:
\par
\medskip
(1.2) Given a pair $(m, n)$ of positive integers, such that $n \mid
m$ and $n$ is divisible by any $p \in \mathbb P$ dividing $m$, there
is a field $F$ and $D \in d(F)$ with ind$(D) = m$ and exp$(D) = n$
(Brauer, see \cite{P}, Sect. 19.6). One can take as $F$ any rational
(i.e. purely transcendental) extension of infinite transcendency
degree over an arbitrary field $F _{0}$.
\par
\medskip
A field $E$ is said to be of Brauer $p$-dimension Brd$_{p}(E) = n$,
where $n \in \mathbb Z$, if $n$ is the least integer for which
ind$(D) \le {\rm exp}(D) ^{n}$ whenever $D \in d(E)$ and $[D] \in
{\rm Br}(E) _{p}$. We say that Brd$_{p}(E) = \infty $, if there
exists a sequence $D _{\nu } \in d(E)$, $\nu \in \mathbb N$, such
that $[D _{\nu }] \in {\rm Br}(E) _{p}$ and ind$(D _{\nu }) > {\rm
exp}(D _{\nu }) ^{\nu }$, for each index $\nu $. By an absolute
Brauer $p$-dimension (abbr, abrd$_{p}(E)$) of $E$, we mean the
supremum sup$\{{\rm Brd}_{p}(R)\colon \ R \in {\rm Fe}(E)\}$. Here
and in the sequel, Fe$(E)$ denotes the set of finite extensions of
$E$ in a separable closure $E _{\rm sep}$. In what follows, we
denote by trd$(F/E)$ the transcendency degree and $I(F/E)$ stands
for the set of intermediate fields of any extension $F/E$.
\par
\medskip
Clearly, ${\rm Brd}_{p}(E) \le {\rm abrd}_{p}(E)$, for every field
$E$ and $p \in \mathbb P$. It is known that ${\rm Brd}_{p}(E) = {\rm
abrd}_{p}(E) = 1$, for every $p \in \mathbb P$, in the following
cases:
\par
\medskip
(1.3) (i) $E$ is a global or local field (by class field theory,
see, e.g., \cite{CF}, Chs. VI and Ch. VII, by Serre and Tate,
respectively);
\par
(ii) $E$ is the function field of an algebraic surface (defined) over
an algebraically closed field $E _{0}$ \cite{Jong}, \cite{Lieb};
\par
(iii) $E$ is the function field of an algebraic curve over a
pseudo algebraically closed field $E _{0}$ with
cd$_{p}(\mathcal{G}_{E _{0}}) > 0$ \cite{Efr1}.
\par
\medskip
By a Brauer dimension and an absolute Brauer dimension of $E$, we
mean the suprema Brd$(E) = {\rm sup}\{{\rm Brd}_{p}(E)\colon \ p \in
\mathbb P\}$ and abrd$(E) = {\rm sup}\{{\rm abrd}_{p}(E)\colon $ $p
\in \mathbb P\}$, respectively. It would be of interest to know
whether the function fields of algebraic varieties over a global,
local or algebraically closed field are of finite absolute Brauer
dimensions. Note also that fields of finite absolute Brauer
$p$-dimensions, for all $p \in \mathbb P$, are singled out for their 
place in research areas like Galois cohomology (cf. \cite{Kah}, 
Sect. 3, \cite{Ch1}, Remark~3.6, and \cite{Ch2}, the end of Section
3 and Corollary~5.7) and the structure theory of their locally
finite-dimensional central division algebras (see \cite{Ch1},
Proposition~1.1 and the paragraph at the bottom of page 2). These
observations draw one's attention to the following open problem:
\par
\medskip
(1.4) Find whether the class of fields $E$ of finite absolute Brauer
$p$-dimensions, for a fixed $p \in \mathbb P$ different from
char$(E)$, is closed under the formation of finitely-generated
extensions.
\par
\medskip
The following result of \cite{Ch1} is used there for proving that
the class of fields $E$ with Brd$(E) < \infty $ is not closed under
taking finitely-generated extensions:
\par
\medskip
\begin{theo}
Let $E$ be a field, $p \in \mathbb P$ and $F/E$ a finitely-generated
extension, such that {\rm trd}$(F/E) = t \ge 1$. Then:
\par
{\rm (i)} {\rm Brd}$_{p}(F) \ge {\rm abrd}_{p}(E) + t - 1$, if {\rm
abrd}$_{p}(E) < \infty $ and $F/E$ is rational;
\par
{\rm (ii)} When {\rm abrd}$_{p}(E) = \infty $, there are $\{D _{n,m} \in
d(F)\colon \ n \in \mathbb N, m = 1, \dots , n\}$ with {\rm exp}$(D
_{n,m}) = p ^{m}$ and {\rm ind}$(D _{n,m}) = p ^{n}$, for each
admissible pair $(n, m)$;
\par
{\rm (iii)} Brd$_{p}(F) = \infty $, provided $p = {\rm char}(E)$ and
the degree $[E\colon E ^{p}]$ is infinite, where $E ^{p} = \{e
^{p}\colon \ e \in E\}$; if {\rm char}$(E) = p$ and $[E\colon E
^{p}] = p ^{\nu } < \infty $, then $\nu + t - 1 \le {\rm
Brd}_{p}(F) < \nu + t$.
\end{theo}

\medskip
Theorem 1.1 is supplemented in \cite{Ch1}, Sect. 3, as follows:

\medskip
(1.5) Given a finitely-generated field extension $F/E$ with
trd$(F/E) = t \ge 1$ and abrd$_{p}(E) < \infty $ when $p$ runs
across some nonempty subset $P \subseteq \mathbb P$, there exists a
finite subset $P(F/E)$ of $P$, such that {\rm Brd}$_{p}(F) \ge {\rm
abrd}_{p}(E) + t - 1$, for each $p \in P \setminus P(F/E)$.

\medskip
It is worth noting that there exist field extensions $F/E$
satisfying the conditions of (1.5), for $P = \mathbb P$, such that
$P(F/E)$ is necessarily nonempty.

\medskip
{\bf Example.} Let $E$ be a real closed field, $F$ the function
field of the Brauer-Severi variety corresponding to the symbol
$E$-algebra $A = A _{-1}(-1,$ $-1; E)$, and $F ^{\prime } = F
\otimes _{E} E(\sqrt{-1})$. By the Artin-Schreier theory (cf.
\cite{L2}, Ch. XI, Theorem~2), then $E(\sqrt{-1}) = E _{\rm sep}$,
whence abrd$_{p}(E) = 0$, for all $p \in \mathbb P \setminus \{2\}$.
Since $-1$ does not lie in the norm group $N(E(\sqrt{-1})/E)$, it
also follows that $A \in d(E)$. Note further that trd$(F/E) = 1$,
$[A \otimes _{E} F] = 0$ in Br$(F)$, and $F ^{\prime }/E(\sqrt{-1})$
is a rational extension (see \cite{Sal2}, Theorem~13.8 and
Corollaries~13.9 and 13.16). In view of Tsen's theorem (cf.
\cite{P}, Sect. 19.4), the noted property of $F ^{\prime }$ ensures
that it is a $C _{1}$-field, so it follows from \cite{S1}, Ch. II,
Proposition~6, that cd$(\mathcal{G}_{F'}) \le 1$. As $A \otimes _{E}
F \cong A _{1}(-1, -1; F)$ over $F$, the equality $[A \otimes _{E}
F] = 0$ implies $F$ is a nonreal field, so it follows from the
Artin-Schreier theory that $\mathcal{G}_{F}$ is a torsion-free
group. Observing finally that $\mathcal{G}_{F'}$ embeds in
$\mathcal{G}_{F}$ as an open subgroup, one obtains from \cite{S1}, 
Ch. I, 4.2, Corollary~3, that cd$(\mathcal{G}_{F}) \le 1$, which 
means that abrd$(F) = 0 < {\rm abrd}_{2}(E) = 1$.

\medskip
Statement (1.1), Theorem 1.1 and basic properties of
finitely-generated field extensions (cf. \cite{L2}, Ch. X) imply the
following:
\par
\medskip
(1.6) If the answer to (1.4) is affirmative, for some $p \in \mathbb
P$, $p \neq {\rm char}(E)$, and each finitely-generated extension
$F/E$ with trd$(F/E) = t \ge 1$, then there exists $c _{t}(p) \in 
\mathbb N$, such that Brd$_{p}(\Phi ) \le c _{t}(p)$
whenever $\Phi /E$ is a finitely-generated extension and trd$(\Phi
/E) < t$ (see also \cite{Ch1}, Proposition~4.6).

\medskip
Theorem 1.1 (i) shows that the solution to \cite{ABGV}, Problem~4.5,
concerning the possibility to find a good definition of a field
dimension dim$(E)$, is negative except, possibly, in the case of
abrd$(E) < \infty $. In addition, it implies that if abrd$(E) <
\infty $ and \cite{ABGV}, Problem~4.5, is solved affirmatively, for
all finitely-generated extensions $F/E$, then the fields $F$ satisfy
much stronger conditions than the one conjectured by (1.6) (see
\cite{Ch1}, (1.5)). As to our next result (for a proof, see
\cite{Ch1}, Proposition~5.8), it indicates that the answer to (1.4)
will be positive, for finitely-generated extensions $F/E$ with
trd$(F/E) \le n$, for some $n \in \mathbb N$, if this is the case in
zero characteristic (see also \cite{Ch1}, Remark~5.9, for an
application of de Jong's theorem \cite{Jong}):

\medskip
(1.7) Let $E$ be a field of characteristic $q > 0$ and $F/E$ a
finitely-generated extension. Then there exists a field $E ^{\prime
}$ with char$(E) = 0$ and a finitely-generated extension $F
^{\prime }/E ^{\prime }$ satisfying the following:
\par
(i) $\mathcal{G}_{E'} \cong \mathcal{G}_{E}$ and trd$(F
^{\prime }/E ^{\prime }) = {\rm trd}(F/E)$;
\par
(ii) {\rm Brd}$_{p}(F ^{\prime }) \ge {\rm Brd}_{p}(F)$, {\rm
abrd}$_{p}(F ^{\prime }) \ge {\rm abrd}_{p}(F)$, {\rm Brd}$_{p}(E
^{\prime }) = {\rm Brd}_{p}(E)$ and {\rm abrd}$_{p}(E ^{\prime }) =
{\rm abrd}_{p}(E)$, for each $p \in \mathbb P$ different from $q$.

\medskip
The proof of Theorem 1.1 in \cite{Ch1} relies on the following two
lemmas. When $\mu = 1$, the former one is a theorem due to Albert. 
Besides in \cite{Ch1}, Sect. 3, a proof of the former lemma can be 
found in \cite{PaSu}, Sect. 1.

\medskip
\begin{lemm}
A field $E$ satisfies the inequality {\rm abrd}$_{p}(E) \le \mu $,
for some $p \in \mathbb P$ and $\mu \in \mathbb N$, if and only if,
for each $E ^{\prime } \in {\rm Fe}(E)$, {\rm ind}$(\Delta _{E'})
\le p ^{\mu }$ whenever $\Delta _{E'} \in d(E ^{\prime })$ and {\rm
exp}$(\Delta ) = p$.
\end{lemm}

\medskip
\begin{lemm}
Let $K$ be a field, $F = K(X)$ a rational extension of $K$ with {\rm
trd}$(F/K) = 1$, $f(X) \in K[X]$ a separable and irreducible
polynomial over $K$, $L$ an extension of $K$ in $K _{\rm sep}$
obtained by adjunction of a root of $f$, $v$ a discrete valuation of
$F$ acting trivially on $K$ with a uniform element $f$, and $(F
_{v}, \bar v)$ a Henselization of $(F, v)$. Suppose that $\widetilde
D \in d(L)$ is an algebra of exponent $p$. Then $L$ is
$K$-isomorphic to the residue field of $(F _{v}, \bar v)$, and there
exist $D ^{\prime } \in d(F _{v})$ and $D \in d(F)$, such that {\rm
exp}$(D) = {\rm exp}(D ^{\prime }) = p$, $[D \otimes _{F} F _{v}] =
[D ^{\prime }]$, and $D ^{\prime }$ is an inertial lift of
$\widetilde D$ over $F _{v}$.
\end{lemm}

\medskip
\section{\bf The main result}

The purpose of this paper is to prove the following assertion which
applied to a field with abrd$_{p}(E) = 0$, improves the inequality in
Theorem 1.1 (i):

\medskip
\begin{theo}
Let $F$ be a finitely-generated extension of a field $E$ with {\rm
cd}$_{p}(\mathcal{G}_{E}) \neq 0$. Then {\rm Brd}$_{p}(F) \ge {\rm
trd}(F/E)$ except, possibly, when $p = 2$, the Sylow
pro-$2$-subgroups of $\mathcal{G}_{E}$ are of order $2$, and $F$ is
a nonreal field.
\end{theo}

\medskip
The following result is contained in \cite{Ch1}, Propositions~4.6
and 5.10, and is obtained by the method of proving Theorem 2.1 (see
also \cite{Ch2}, (4.10) and Proposition~4.3):

\medskip
\begin{theo}
Assume that $E$ is a field of type pointed out in (1.3). Then {\rm
Brd}$_{p}(F) \ge 1 + {\rm trd}(F/E)$, for every finitely-generated
extension $F/E$.
\end{theo}

\medskip
\begin{rema}
(i): Theorem 2.1 ensures that Brd$_{p}(\Phi ) \ge n$, $p \in \mathbb
P$, if $\Phi $ is a finitely-generated extension of a quasifinite 
field $\Phi _{0}$, and trd$(\Phi /\Phi _{0}) = n$. Therefore, one 
obtains following the proof of \cite{Ch1}, Proposition~5.10, that 
the conclusion of Theorem 2.2 remains valid, if $E$ is endowed with 
a Henselian discrete valuation whose residue field is quasifinite. 
\par
(ii): Given a finitely-generated field extension $F/E$ with
trd$(F/E) = k$, Theorem 2.1 implies Nakayama's inequalities 
Brd$_{p}(F) \ge k - 1$, $p \in \mathbb P$ (cf. \cite{Jong}, Sect. 
2). When cd$_{p}(\mathcal{G}_{E}) = 0$, for some $p$, and $E$ is 
perfect in case $p = {\rm char}(E)$, we have Brd$_{p}(F) = k - 1$ if 
and only if this holds in case $E$ is algebraically closed. The 
claim that Brd$(F) = k - 1$ when $E$ is algebraically closed is the 
content of the so called Standard Conjecture, for function fields of 
algebraic varieties over an algebraically closed field (see 
\cite{Lieb}, Sect. 1, \cite{LiKr}, page 3, and for relations with 
(1.4), the end of \cite{Ch1}, Sect. 4).
\end{rema}

\medskip
The proof of Theorem 2.1 is based on the same idea as the one of
Theorem 1.1. It relies on the following lemmas proved in \cite{Ch1}.
\par
\medskip
\begin{lemm}
Let $(K, v)$ be a nontrivially real-valued field, and $(K _{v}, \bar
v)$ a Henselization of $(K, v)$. Assume that $\Delta _{v} \in d(K
_{v})$ has exponent $p \in \mathbb P$. Then there exists $\Delta \in
d(K)$, such that {\rm exp}$(\Delta ) = p$ and $[\Delta \otimes _{K}
K _{v}] = [\Delta _{v}]$.
\end{lemm}

\medskip
Lemma 2.4 is essentially due to Saltman \cite{Sal1}, and our next
lemma is a special case of the Grunwald-Wang theorem (cf. \cite{LR},
Theorems~1 and 2).

\medskip
\begin{lemm}
Let $F$ be a field, $S = \{v _{1}, \dots , v _{s}\}$ a finite set of
non-equivalent nontrivial real-valued valuations of $F$, and
for each index $j$, let $F _{v _{j}}$ be a Henselization of $K$ in $K
_{\rm sep}$ relative to $v _{j}$, and $L _{j}/F _{v _{j}}$ be a cyclic
field extension of degree $p ^{\mu _{j}}$, for some $p \in P$ and $\mu
_{j} \in \mathbb N$. Put $\mu = {\rm max}\{\mu _{1}, \dots , \mu _{s}\}$
and suppose that $\sqrt{-1} \in F$ in case $\mu \ge 3$, $p = 2$ and
{\rm char}$(F) = 0$. Then there exists a degree $p ^{\mu }$ cyclic field
extension $L/F$, such that $L _{v _{j}'}$ is $F _{v
_{j}}$-isomorphic to $L _{j}$, where $v _{j} ^{\prime }$ is a valuation
of  $L$ extending $v _{j}$, for $j = 1, \dots , s$.
\end{lemm}

\medskip
In the rest of this Section, we recall some general results on
Henselian valuations which are used (often implicitly, like Lemma
1.3) for proving Theorem 2.1. A Krull valuation $v$ of a field $K$
is called Henselian, if $v$ extends uniquely, up-to an equivalence,
to a valuation $v _{L}$ on each algebraic extension $L$ of $K$.
Assuming that $v$ is Henselian, denote by $v(L)$ the value group and
by $\widehat L$ the residue field of $(L, v _{L})$. It is known that
$\widehat L/\widehat K$ is an algebraic extension and $v(K)$ is a
subgroup of $v(L)$. When $L/K$ is finite and $e(L/K)$ is the index 
of $v(K)$ in $v(L)$, by Ostrowski's theorem \cite{Efr2}, 
Theorem~17.2.1, $[\widehat L\colon \widehat K]e(L/K)$ divides 
$[L\colon K]$ and $[L\colon K][\widehat L\colon \widehat K] 
^{-1}e(L/K) ^{-1}$ is not divisible by any $p \in \mathbb P$, $p 
\neq {\rm char}(\widehat K)$. In particular, if char$(\widehat K)$ 
does not divide $[L\colon K]$, then $[L\colon K] = [\widehat L\colon 
\widehat K]e(L/K)$. Ostrowski's theorem implies that there are group 
isomorphisms $v(K)/pv(K) \cong v(L)/pv(L)$, $p \in \mathbb P$, and 
in case char$(\widehat K) \dagger [L\colon K]$, they are canonically 
induced by the natural embedding of $K$ into $L$.
\par
\medskip
As usual, a finite extension $R$ of $K$ is called inertial, if $[R\colon
K] = [\widehat R\colon \widehat K]$ and $\widehat R$ is separable
over $\widehat K$. It follows from the Henselity of $v$ that the
compositum $K _{\rm ur}$ of inertial extensions of $K$ in $K _{\rm
sep}$ has the following properties:
\par
\medskip
(2.1) (i) $v(K _{\rm ur}) = v(K)$ and finite extensions of $K$ in $K
_{\rm ur}$ are inertial;
\par
(ii) Each finite extension of $\widehat K$ in $\widehat K _{\rm
sep}$ is $\widehat K$-isomorphic to the residue field of an inertial
extension of $K$; hence, $\widehat K _{\rm ur}$ is $\widehat
K$-isomorphic to $\widehat K _{\rm sep}$;
\par
(iii) $K _{\rm ur}/K$ is a Galois extension with $\mathcal{G}(K
_{\rm ur}/K) \cong \mathcal{G}_{\widehat K}$.

\medskip
Similarly, it is known that each $\Delta \in d(K)$ has a unique,
up-to an equivalence, valuation $v _{\Delta }$ extending $v$ so that
the value group $v(\Delta )$ of $(\Delta , v _{\Delta })$ is abelian
(see \cite{JW}). Note that $v(\Delta )$ includes $v(K)$ as an
ordered subgroup of index $e(\Delta /K) \le [\Delta \colon K]$, the
residue division ring $\widehat {\Delta }$ of $(\Delta , v _{\Delta
})$ is a $\widehat K$-algebra, and $[\widehat {\Delta }\colon
\widehat K] \le [\Delta \colon K]$. Moreover, by Ostrowski-Draxl's
theorem (cf. \cite{JW}, (1.2)), $e(\Delta /K)[\widehat {\Delta
}\colon \widehat K] \mid [\Delta \colon K]$, and in case
char$(\widehat K) \dagger [\Delta \colon K]$, $[\Delta \colon K] =
e(\Delta /K)[\widehat \Delta \colon \widehat K]$. An algebra $D \in
d(K)$ is called inertial, if $[D\colon K] = [\widehat D\colon
\widehat K]$ and $\widehat D \in d(\widehat K)$. In what follows, we
also need the following results (see \cite{JW}, Remark~3.4 and
Theorems~2.8 and 3.1):
\par
\medskip
(2.2) (i) Each $\widetilde D \in d(\widehat K)$ has a unique, up-to
an $F$-isomorphism, inertial lift $D$ over $K$(i.e. $D \in d(K)$,
$D$ is inertial over $K$ and $\widehat D = \widetilde D$).
\par
(ii) The set IBr$(K)$ of Brauer equivalence classes of inertial
$K$-algebras forms a subgroup of Br$(K)$ canonically isomorphic
to Br$(\widehat K)$.
\par
(iii) For each $\Theta \in d(K)$ inertial over $K$, and any $R \in
I(K _{\rm ur}/K)$, $[\Theta \otimes _{K} R] \in {\rm IBr}(R)$ and
ind$(\Theta \otimes _{K} R) = {\rm ind}(\widehat \Theta \otimes
_{\widehat K} \widehat R)$.

\medskip
\section{\bf Proof of Theorem 2.1}

\medskip
Let $E$ be a field with cd$_{p}(\mathcal{G}_{E}) > 0$, for some $p
\in \mathbb P$, and let $F/E$ be a finitely-generated extension.
Throughout this Section, $E _{\rm sep}$ is identified with its
$E$-isomorphic copy in $F _{\rm sep}$, and for any field $Y$, $r
_{p}(Y)$ denotes the rank of the Galois group $\mathcal{G}(Y(p)/Y)$
of the maximal $p$-extension $Y(p)$ of $Y$ (in $Y _{\rm sep}$) as a
pro-$p$-group. Assuming that trd$(F/E) = t$ and $G _{p}$ is a
Sylow pro-$p$-subgroup of $\mathcal{G}_{E}$, we deduce Theorem 2.1
by proving the following:
\par
\medskip
(3.1) There exists $D \in d(F)$, such that exp$(D) = p$, ind$(D) = p
^{t}$ and $D$ is presentable as a tensor product of cyclic
$F$-algebras of degree $p$ except, possibly, in the case where $p =
2$, $G _{2}$ is of order $2$ and $F$ is a nonreal field.
\par
\medskip
Let $E _{p}$ be the fixed field and $o(G _{p})$ the order of $G
_{p}$. Our assumptions show that $r _{p}(E _{p}) \ge 1$, which
implies the existence of a field $M \in {\rm Fe}(E)$ with $r _{p}(M)
\ge 1$ (apply the method of proving \cite{P}, Sect. 13.2,
Proposition~b). Moreover, $M$ can be chosen to be nonreal unless $p
= 2$ and $o(G _{p}) = 2$. Assuming that $M$ is nonreal, one obtains
from \cite{Wh}, Theorem~2, that there exists a $\mathbb Z
_{p}$-extension $\Phi $ of $M$ in $E _{\rm sep}$. Hence, by Galois
theory and the fact that $\mathbb Z _{p}$ is continuously isomorphic
to its open subgroups, $\Phi M ^{\prime }/M ^{\prime }$ is Galois
with $\mathcal{G}(\Phi M ^{\prime }/M ^{\prime }) \cong \mathbb Z
_{p}$, for each $M ^{\prime } \in {\rm Fe}(E)$. This makes it easy 
to obtain from basic properties of valuation prolongations on finite 
extensions that $M$ can be chosen as an $E$-isomorphic copy of the
residue field of a height $t$ valuation $v$ of $F$, trivial on
$E$ with $v(F) = \mathbb Z ^{t}$. Here $\mathbb Z ^{t}$ is viewed as
an ordered abelian group with respect to the inversely-lexicographic
ordering.
\par
Let $(F _{v}, \bar v)$ be a Henselization of $(F, v)$. Suppose first
that $t = 1$ and take $\pi \in F$ so that $\langle v(\pi
)\rangle = v(F)$. Then $v$ lies in an infinite system of
nonequivalent discrete valuations of $F$ trivial on $E$ (cf.
\cite{CF}, Ch. II, Lemma~3.1). In view of Lemma 2.5, this implies
the existence of degree $p$ cyclic extensions $F _{n}$, $n \in
\mathbb N$, of $F$, such that $F _{1}/F$ is inertial relative to
$v$, and $F _{n} \subset F _{v}$, $n \ge 2$. Let $\varphi _{n}$ be a
generator of $\mathcal{G}(F _{n}/F)$, for each $n \in \mathbb N$. It
follows from the choice of $F _{1}$ that the cyclic $F$-algebra $(F
_{1}/F, \sigma _{1}, \pi )$ lies in $d(F)$ and $(F _{1}/F, \sigma
_{1}, \pi ) \otimes _{F} F _{v} \in d(F _{v})$, which proves (3.1)
in case $t = 1$.
\par
Assume now that $t \ge 2$, and fix elements $\pi _{1}, \dots ,
\pi _{t} \in K$ so that $v(F)$ be generated by the set
$\{v(\pi _{j})\colon \ j = 1, \dots , t\}$, and $H = \langle
v(\pi _{1})\rangle $ be the minimal nontrivial isolated subgroup of
$v(F)$. Then $v$ and $H$ induce canonically on $F$ a valuation $v
_{H}$ with $v _{H}(F) = v(F)/H$; also, they give rise to a valuation
$\hat v _{H}$ of the residue field $F _{H}$ of $(F, v _{H})$ with
$\hat v _{H}(F _{H}) = H$ and a residue field equal to $M$ (cf.
\cite{Efr2}, Sect. 5.2). In addition, it is easily verified that $F
_{H}/E$ is a finitely-generated extension with trd$(F _{H}/E) = 1$.
Hence, by the proof of the already considered special case of
Theorem 2.1, there exist $D _{H} \in d(F _{H})$ and $\Psi _{H} \in
I(F _{H,{\rm sep}}/F _{H})$, such that ind$(D _{H}) = p$, $D _{H}
\otimes _{F _{H}} \Psi _{H} \in d(\Psi _{H})$, and $I(\Psi _{H}/F
_{H})$ contains infinitely many degree $p$ cyclic extensions of $F
_{H}$. Observing now that $v _{H}$ is of height $t - 1$, and
using repeatedly (2.2), Lemmas 2.4, 2.5 and \cite{Ch2}, (3.1) (i),
one proves that there exists a cyclic $F$-algebra $D ^{\prime } \in
d(F)$, such that ind$(D ^{\prime }) = p$, $D ^{\prime } \otimes _{F}
F _{v _{H}} \in d(F _{v _{H}})$ and $D ^{\prime } \otimes _{F} F _{v
_{H}}$ is an inertial lift of $D _{H}$ over a Henselization $F _{v
_{H}}$ of $F$ relative to $v _{H}$. Similarly, it can be deduced
from (2.2) that each degree $p$ cyclic extension of $F _{H}$ is
realizable as the residue field of an inertial cyclic degree $p$
extension of $F$ relative to $v _{H}$. This implies the existence of
an inertial extension $(F ^{\prime }, v _{H} ^{\prime })/(F, v
_{H})$, such that $[F ^{\prime }\colon F] = [F ^{\prime }F _{v
_{H}}\colon F _{v _{H}}] = p ^{t-1}$, $D ^{\prime } \otimes _{F} F
^{\prime } \in d(F)$ and $F ^{\prime } = F _{2} \dots F _{t}$,
where $F _{i}/F$ is a degree $p$ cyclic extension of $F$, for $i =
2, \dots , t$. In view of Morandi's theorem (cf. \cite{JW},
Proposition~1.4), it is now easy to construct an algebra $\Delta \in
d(F)$, such that exp$(\Delta ) = p$, ind$(\Delta ) = p ^{t-1}$,
$\Delta $ is presentable as a tensor product of cyclic $F$-algebras
of degree $p$, $\Delta \otimes _{F} F _{v _{H}} \in d(F _{v _{H}})$,
$\Delta \otimes _{F} F _{v _{H}}$ is nicely semi-ramified over $F
_{v _{H}}$, in the sense of \cite{JW}, and $(D ^{\prime } \otimes
_{F} \Delta ) \otimes _{F} F _{v _{H}} \in d(F _{v _{H}})$.
Therefore, $D ^{\prime } \otimes _{F} \Delta \in d(F)$, exp$(D
^{\prime } \otimes _{F} \Delta ) = p$ and ind$(D ^{\prime } \otimes
_{F} \Delta ) = p ^{t}$, which proves (3.1), under the hypothesis
that $o(G _{p}) > 2$.
\par
It remains to be seen that (3.1) holds when $p = 2$, $F$ is formally
real and $o(G _{2}) = 2$. By the Artin-Schreier theory, $o(G _{2}) =
2$ if and only if the fixed field $E _{2}$ is real closed. Our proof
also relies on the following lemma.

\medskip
\begin{lemm}
Let $E$ be a formally real field and $F$ a finitely-generated
extension of $E$ with {\rm trd}$(F/E) = 1$. Then $F$ is formally
real if and only if it has a discrete valuation $v$ trivial on $E$,
whose residue field $\widehat F$ is formally real.
\end{lemm}

\medskip
\begin{proof}
It is known and easy to prove (cf. \cite{L1}, Lemma~1) that if $F$
is a nonreal field and $\omega $ is a discrete valuation of $F$
trivial on $E$, then the residue field of $(F, \omega )$ is nonreal
as well. Assume now that $F$ is formally real, fix a real closure $F
^{\prime }$ of $F$ in $F _{\rm sep}$, and put $E ^{\prime } = E
_{\rm sep} \cap F ^{\prime }$. Observe that $E _{\rm sep}F ^{\prime
}/F ^{\prime }$ is a Galois extension with $\mathcal{G}(E _{\rm
sep}F ^{\prime }/F ^{\prime }) \cong \mathcal{G}_{E'}$. Since, by
the Artin-Schreier theory, $F _{\rm sep} = F ^{\prime }(\sqrt{-1}) =
E _{\rm sep}F ^{\prime }$, this means that $E _{\rm sep} = E
^{\prime }(\sqrt{-1})$, whence, $E ^{\prime }$ is a real closure of
$E$ in $E _{\rm sep}$. Note also that $E ^{\prime }F/E ^{\prime }$
is finitely-generated, trd$(E ^{\prime }F/E ^{\prime }) = 1$ and $E
^{\prime }F \subseteq F ^{\prime }$, i.e. the extension $E ^{\prime
}F/E ^{\prime }$ satisfies the conditions of Lemma 3.1. This enables
one to deduce from \cite{L1}, Theorem~6 and Proposition, that $E
^{\prime }F$ has a discrete valuation $v ^{\prime }$ trivial on $E
^{\prime }$ and with a residue field $E ^{\prime }$. It is now easy
to see that the valuation $v$ of $F$ induced by $v ^{\prime }$ has
the properties required by Lemma 3.1. Specifically, $\widehat F$ is
$E$-isomorphic to a finite extension of $E$ in $E ^{\prime }$.
\end{proof}

\medskip
We are now in a position to prove the remaining case of (3.1).
Suppose first that $t = 1$, put $F _{0} = E(X)$, for some $X
\in F$ transcendental over $E$, and denote by $\Omega _{0}$ the
extension of $F _{0}$ in $F _{\rm sep}$ generated by the square
roots of the totally positive elements of $F _{0}$ (i.e. those
realizable over $F _{0}$ as finite sums of squares, see \cite{L2},
Ch. XI, Proposition~2). Then $F\Omega _{0}$ is formally real,
which implies $A _{-1}(-1, -1; \Omega ) \in d(\Omega )$, for each
$\Omega \in I(F\Omega _{0}/F _{0})$, proving the assertion of (3.1).
Note also that $\Omega _{0}/F _{0}$ is an infinite Galois extension
with $\mathcal{G}(\Omega _{0}/F _{0})$ of exponent $2$. This follows
from Kummer theory and the fact that cosets $(X ^{2} + a ^{2})F _{0}
^{\ast 2}$, $a \in E ^{\ast }$, generate an infinite subgroup of $F
_{0} ^{\ast }/F _{0} ^{\ast 2}$.
\par
Assume now that $t \ge 2$, define $F _{0}$ and $\Omega _{0}$ as
above and denote by $F _{1}$ the algebraic closure of $F _{0}$ in
$F$. Applying Lemma 3.1 and proceeding by induction on $t $, one
concludes that $F$ has valuation $v$ trivial on $F _{1}$, such that
$v(F) = \mathbb Z ^{t-1}$, $v$ is of height $t - 1$ and $\widehat F$
is a formally real finite extension of $F _{1}$. Fix a Henselization
$(F _{v}, \bar v)$ and an $F _{1}$-isomorphic copy $F _{1} ^{\prime
}$ of $\widehat F$ in $F _{\rm sep}$. It is easily verified that $F
_{1} ^{\prime }\Omega _{0}$ is a formally real field and $F\Omega
_{0}/F$ is a Galois extension. As $v$ is of height $t - 1$, one
proves, using repeatedly (2.1) and Lemma 2.5, that $I(F\Omega
_{0}/F)$ contains infinitely many quadratic and inertial extensions
of $F$ relative to $v$. Therefore, there exist fields $Y _{n} \in
I(F\Omega _{0}/F)$, $n \in \mathbb N$, such that $[Y _{n}\colon F] =
2$, $[Y _{1} \dots Y _{n}\colon F] = 2 ^{n}$ and $Y _{1} \dots Y
_{n}$ is inertial over $F$ relative to $v$, for each index $n$. Fix
a generator $q _{j}$ of $\mathcal{G}(Y _{j}/F)$, and take elements
$\pi _{j} \in F$, $j = 2, \dots , t$, so that $\langle v(\pi _{2}),
\dots , v(\pi _{t})\rangle = v(F)$. Put $\Delta _{1} = A _{-1}(-1,
-1; F)$ and consider the cyclic $F$-algebras $\Delta _{j} = (Y
_{j}/F, q _{j}, \pi _{j})$, $j = 2, \dots , t$. Since $F _{1}
^{\prime }\Omega _{0}$ is formally real, $A _{-1}(-1, -1; F) \otimes
_{F} F _{1} ^{\prime }\Omega _{0} \in d(F _{1} ^{\prime }\Omega
_{0})$, so it follows from Morandi's theorem, the noted properties
of the fields $Y _{n}$, $n \in \mathbb N$, and the choice of $\pi
_{2}, \dots , \pi _{t}$, that the $F$-algebra $\Delta = \Delta _{1}
\otimes _{F} \dots \otimes _{F} \Delta _{t}$ lies in $d(F)$ and
$\Delta \otimes _{F} F _{v} \in d(F _{v})$. This yields exp$(\Delta
) = 2$ and ind$(\Delta ) = 2 ^{t}$, so (3.1) and Theorem 2.1 are
proved.

\par

\end{document}